\documentclass[10pt]{article}
\usepackage{graphicx}
\usepackage{epsfig,amssymb,color}

\begin{document}
\baselineskip=15pt
\DeclareGraphicsExtensions{.pdf,.png,.jpg}
\begin{center}

\def\d{{\mathrm{d}}}

{\huge Singularity in the discrete-time model of 
  impacting mechanical systems}

\vskip .5cm

{\bf Soumya Kundu,
\footnote{Presenting author: 
Phone: +91-9433217289, Fax: +91-3222282262,
email: soumya.joy@gmail.com} and 
Soumitro Banerjee}

{\it Department of Electrical Engineering, Indian Institute of
  Technology Kharagpur, India}

\end{center}

\begin{abstract}
It is known that many peculiar nonlinear vibration problems in
impacting systems are caused by grazing incidences. Such bifurcation
phenomena are normally investigated through the Poincar\'e map. The
discrete-time map of a simple impact oscillator was derived by
Nordmark, which showed that there should be a square-root singularity
in the Jacobian matrix close to the grazing condition. In this paper
we show that the square root singularity will be expressed only in the
trace of the Jacobian matrix, while the determinant remains invariant
across the grazing condition.
\end{abstract}


\vspace{0.1in}

\section{Introduction}

Mechanical systems with impacts between elements occur frequently in
engineering practice. Vibration problems in such systems essentially
hinge on the dynamics of a moving body, possibly in a
mass-spring-damper combination, impacting with a rigid stop. It
is known that much of the dynamical phenomena in such systems stem
from the conditions pertaining to the situation when one body just
grazes the other. That is why much attention has been given to the
grazing condition, and the bifurcation phenomena resulting from that
\cite{shaw-holmes,Peterka-Vacik,Kobrynskii,nordmark}.

It is convenient to analyze bifurcation phenomena in any dynamical
system by obtaining a discrete time model or map by the method of
Poincar\'e surface of section. The structure of the obtained map
determines the dynamics of the physical system. Thus, in understanding
the dynamics of impacting systems, researchers have tried to obtain
the structure of the map. It is obvious that the map should be
piecewise smooth, since the map for non-impacting condition and that
for the impacting condition should be different, and the two should be
separated by the grazing condition. Most important in the respect is
the question: What is the structure of the map in the neighborhood of
a grazing orbit?

Nordmark \cite{nordmark} first addressed this question. He showed that
in the non-impacting side the map is linear, while in the impacting
side it has a square-root term. This implies that the derivative of
the map approaches infinity as the grazing condition is approached
from the impacting side. This results in an infinite stretching of the
state space---which has come to be known as ``square-root
singularity.'' Many researchers studied the behavior of the impact
oscillator based on this map \cite{Chin95,chin-grebogi94,budd1,budd2}.

In a two-dimensional oscillator (mass-spring-damper type), a sampling
in synchronism with the external forcing finction yields a
two-dimensional map. Nordmark derived the condition on the whole
Jacobian matrix, but not on the elements of the matrix. In the
meantime, many other switching dynamical systems---most notably the
power electronic circuits---were found to exhibit a new class of
bifurcation that occurs when a fixed point crosses the border between
two smooth regions in a piecewise smooth map. The development of the
theory of such border collision bifurcation based on a normal form
\cite{pre2d}. The normal form is expressed as
\begin{eqnarray}
\left(\begin{array}{c} x_{k+1} \\ y_{k+1} \end{array} \right) &=&
 \left\{\begin{array}{cc}
\underbrace{\left(\begin{array}{cc}  \tau_L & 1 \\ -\delta_L & 0
\end{array} \right)}_{{\bf J}_L}
 \left(\begin{array}{c} x_k\\ y_k \end{array} \right)
+ \left(\begin{array}{c} 1 \\ 0 \end{array} \right) \mu,& x_k\leq 0 \\
\underbrace{\left(\begin{array}{cc}  \tau_R & 1 \\ -\delta_R & 0
\end{array} \right)}_{{\bf J}_R} \left(\begin{array}{c} x_k \\ y_k
\end{array} \right) + \left(\begin{array}{c} 1 \\ 0 \end{array}
\right) \mu, & x_k\geq 0
\end{array} \right.
\label{eq:2DNF}
\end{eqnarray}
where $\tau_L$ is the trace and $\delta_L$ is the determinant of
the Jacobian matrix ${\bf J}_L$ of the system at a fixed point in one
side and
$\tau_R$ is the trace and $\delta_R$ is the determinant of the
Jacobian matrix ${\bf J}_R$ of the system evaluated at a fixed
point in the other side.

The natural question in relation
with the impact oscillator was: How does the trace and the determinant
change as an impact oscillator is driven from a non-impacting state to
an impacting state with the change of a parameter?

In this paper we probe this issue, and
analytically prove that the determinant must be invariant while the
trace alone should exhibit the square-root singularity.

\section{Impacting Hybrid System Description}

\begin{figure}[tbh]
\begin{center}
\includegraphics[width=2in]{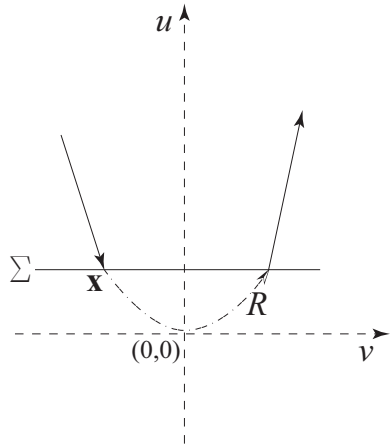}
\caption{Upon impact, the velocity instantly reverses while the
  position remains the same. Thus the state instantly jumps to a new
  position.} \label{impact1}
\end{center}
\end{figure}

An impacting system (Fig.~\ref{impact1}) is governed by a set of
ordinary differential equations (ODEs) coupled with a set of reset maps as 

\begin{displaymath}
\dot{\bf x}=F({\bf x}), \;\;\; \mbox{if}\;\;\; {\bf x}\in S^+ 
\end{displaymath}
\begin{equation}
\label{eqn1}
{\bf x} \mapsto R({\bf x}), \;\;\; \mbox{if}\;\;\;{\bf x}\in\Sigma
\end{equation}
where, $S^+={\{}{\bf x}:H({\bf x})>0{\}}$ and $\Sigma={\{}{\bf
  x}:H({\bf x})=0{\}}.$ $H({\bf x})$ is a smooth function, zero set of
which defines the hard boundary $\Sigma$. The flow given by
(\ref{eqn1}) is restricted only in the region
$S^+\cup\Sigma$.

\noindent

\noindent
Let us now define the normal velocity $v({\bf x})$ as the rate at which the trajectory approaches the impact boundary. It is given by 
\begin{displaymath}
v({\bf x}):=\frac{dH}{dt}=\frac{\partial H}{\partial x} \frac{dx}{dt}  = H_{\bf x}F.
\end{displaymath}

\noindent
Similarly the normal acceleration $a({\bf x})$ of the flow with respect to the boundary is 
\begin{displaymath}
a({\bf x}):=(H_{\bf x}F)_{\bf x}F.
\end{displaymath}

\noindent
We may now be more specific about the form the reset map $R({\bf x})$ takes. To that end, we observe that the reset map has to be a smooth function of the normal velocity $v({\bf x})$ and furthermore $R$ maps to itself when grazing occurs. Since at grazing the normal velocity with respect to the boundary becomes zero ($v({\bf x})=0$), the reset map can be formulated as 
\begin{equation}
\label{eqn2}
R({\bf x})={\bf x}+W({\bf x})v({\bf x})
\end{equation}
where $W$ is a smooth $2\times 1$ matrix.

\section{Grazing and Discontinuity Mapping}

\begin{figure}[tbh]
\begin{center}
\includegraphics[width=2in]{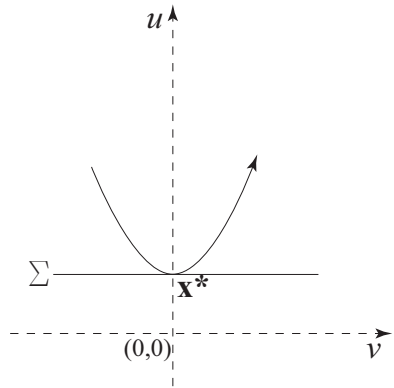}
\caption{Grazing of discontinuity boundary.} \label{impact2}
\end{center}
\end{figure}

Grazing occurs when a trajectory becomes tangent to the discontinuity
boundary $\Sigma$, as shown in Fig.~\ref{impact2}. A point ${\bf x}={\bf x^*}$
is called a regular grazing point if it satisfies the conditions
\begin{displaymath}H({\bf x^*})=0\end{displaymath}
\begin{displaymath}v({\bf x^*})=0\end{displaymath}
\begin{displaymath}a({\bf x^*})=a^*>0\end{displaymath}

\noindent
In addition the scalar function $H({\bf x})$ is assumed to be well
defined at ${\bf x}={\bf x^*}$, i.e., $H_{\bf x}({\bf x^*})\neq0$.

For the part of the flow that does not have any impact with the
discontinuity boundary the mapping is given by the ODE only. Whenever
there is an impact with the boundary, the reset map comes into action
and there is a discontinuity in the flow. The discontinuity near
grazing is of particular interest. Special kinds of mapping have been
proposed to account for this discontinuity \cite{mariobook}. In this
present case the zero-time discontinuity mapping (ZDM) is dealt with.

\begin{figure}[tbh]
\begin{center}
\includegraphics[width=3in]{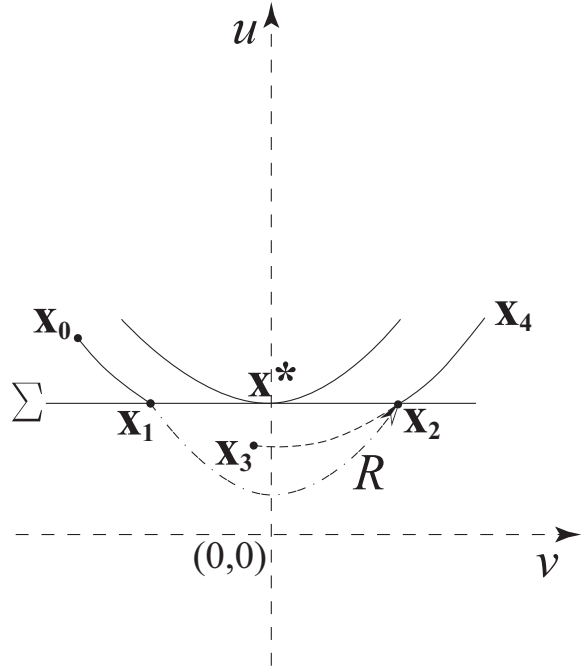}
\caption{Zero-time discontinuity mapping near grazing.} \label{impact3}
\end{center}
\end{figure}

\noindent
Let us consider the situation as shown in Fig.~\ref{impact3}. There is
an orbit which grazes the discontinuity boundary $\Sigma$ at a point
${\bf x^*}$ at some point of time $t_0$. Now let there be another
trajectory (${\bf x_0}{\bf x_1}{\bf x_2}{\bf x_4}$) close to the
grazing orbit. Let us back-trace the trajectory, governed by the ODE
as in (\ref{eqn1}), from the point ${\bf x_2}$ to the point ${\bf
  x_3}$ such that the time taken by the trajectory to reach from ${\bf
  x_0}$ to ${\bf x_2}$ is the same as would have been taken by the
flow to reach from ${\bf x_3}$ to ${\bf x_2}$. Thus, we can consider
the systems' dynamics as if the switching boundary were not there. In
that case we have to assume an instantaneous jump of the state from
${\bf x_0}$ to ${\bf x_3}$.  The ZDM is defined as the mapping ${\bf
  x_0}\mapsto{\bf x_3}$.

Let $\phi({\bf
  x_0},t)$ be the flow obtained as the solution of the ODE starting
from ${\bf x_0}$, i.e. $\phi({\bf x_0},0)={\bf x_0}$. It has been shown in \cite{mariobook} that the form of the ZDM, excluding
higher order terms, is
\begin{equation}
\label{eqn3}
{\bf x_3}={\bf x_0}-W^*({\sqrt{2a^*}})y
\end{equation}
where, $W^*=W({\bf x^*})$, and $y={\sqrt{-H_{\rm min}({\bf
      x_0})}}$.  $H_{\rm min}({\bf x_0})$ is defined as the minimum value
of $H({\phi}({\bf x_0},t))$ with the smallest $|t|$, i.e., the lowest
point that the trajectory would have reached if the switching boundary
were not there.  Obviously for
the situation as described in Fig.~\ref{impact3}, $H_{\rm min}({\bf x_0})$ will be
negative except for the case when ${\bf x_0}$ is the same as ${\bf
  x^*}$.

\begin{figure}[tbh]
\begin{center}
\includegraphics[width=2in]{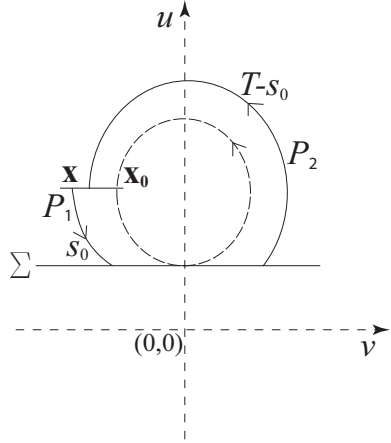}
\caption{A grazing periodic orbit and stroboscopic map.} \label{impact4}
\end{center}
\end{figure}

Now, let
us consider a periodic orbit which has an intersection with the
discontinuity boundary very close to the grazing orbit, as shown in
Fig.~\ref{impact4}. The stroboscopic Poincar\'e map in this case is
 $P_s=P_2{\circ}ZDM{\circ}P_1$, where $P_1$ is the
map that takes a point on the Poincar\'e plane and maps it to
the discontinuity boundary $\Sigma$ by evolution through the ODE in
(\ref{eqn1}), and $P_2$ is the map that takes a point on the
discontinuity boundary $\Sigma$ and maps it back to the
Poincar\'e plane via the ODE. The form of this stroboscopic map can be
derived, in first order approximation, as:
\begin{eqnarray}
P_1 & : & {\bf x}\;{\mapsto}\;N_1{\bf x} \nonumber \\
ZDM{\circ}P_1 & : & {\bf x}\;{\mapsto}\;N_1{\bf x}-{\sqrt{2a^*}}{\sqrt{-H_{\rm min}(N_1{\bf x})}}W^* \nonumber \\
P_2{\circ}ZDM{\circ}P_1 & : & {\bf x}\;{\mapsto}\;N_2N_1{\bf x}-{\sqrt{2a^*}}{\sqrt{-H_{\rm min}(N_1{\bf x})}}N_2W^* \nonumber \\
\end{eqnarray} 
\noindent
where, $N_1:={\frac{dP_1}{d{\bf x}}{\arrowvert}_{{\bf x}={\bf x_0}}}$,
$N_2:={\frac{dP_2}{d{\bf x}}{\arrowvert}_{{\bf x}={\bf x_0}}}$. 

On the Poincar\'e plane, we define $x_0$ as the origin, and the
stroboscopic map in this case is defined as the map which takes the
the initial deviation on the Poincar\'e plane, $({\bf x} - {\bf
  x_0})$, and maps this deviation again on the Poincar\'e plane. Since
${\bf x}_0$ maps to the point $N_2 N_1 {\bf x}_0$, this
map takes the form
\begin{eqnarray}
\label{eqn4}
({\bf x}-{\bf x_0})\;{\mapsto}\;N_2N_1({\bf x}-{\bf x_0})-{\sqrt{2a^*}}{\sqrt{-H_{\rm min}(N_1{\bf x})}}N_2W^*
\end{eqnarray}
Also
$H_{\rm min}(N_1{\bf x})$ is linearized about grazing to the form 
\begin{displaymath}
H_{\rm min}(N_1{\bf x})=H_{\bf x}N_1({\bf x}-{\bf x_0})+H.O.T
\end{displaymath}

\noindent
As can be seen from the expression of the stroboscopic map, and also of ZDM alone, a square-root term $\sqrt{-H_{\rm min}}$ is present which accounts for the square-root singularity when the Jacobian of the map is considered. The next section deals with this Jacobian and the square-root singularity therein.

\section{Investigating the Trace and the Determinant of the Jacobian for the Singularity}

The Jacobian of the stroboscopic map near grazing would be
\begin{equation}
\label{eqn5}
J=N_2{\circ}J_{ZDM}{\circ}N_1
\end{equation}
\noindent
where $J_{ZDM}$ is the Jacobian of the ZDM given by
\begin{equation}
\label{eqn6}
J_{ZDM}=\frac{{\partial}f}{{\partial}{\bf x}}=I+{\sqrt{2a^*}}.{\frac{W^*H_{\bf x}}{2\sqrt{-H_{\rm min}}}}.
\end{equation}

To arrive at the particular forms $W^*$ and $H_{\bf x}$ would
take, let us concentrate on the one degree-of-freedom impact
oscillator (shown in figure 5). The mass (assumed to be unity without
any loass of generality) is tied with a spring-damper element and is
acted upon by an external force $g(t)$. At a distance $\sigma$ from
the mass there is an impacting wall so that for $u<\sigma$ the motion
of the mass is governed by second-order differential equation
\begin{equation}
\label{eqn7}
{\frac{d^2u}{dt^2}}+2\zeta{\omega}_n\frac{du}{dt}+{{\omega}_n}^2u=g(t),\;
\mbox{for}\;u<\sigma
\end{equation}
\noindent
and at $u=\sigma$ the reset map $R$ is applied. This system is a two-dimensional system, i.e. ${\bf x}\in{{\mathbb{R}}^2}$. Considering the velocity of motion of the mass as $v=\frac{du}{dt}$, the state vector can be written as ${\bf x}=(u,v)^T$.

\begin{figure}[tbh]
\begin{center}
\includegraphics[width=2.5in]{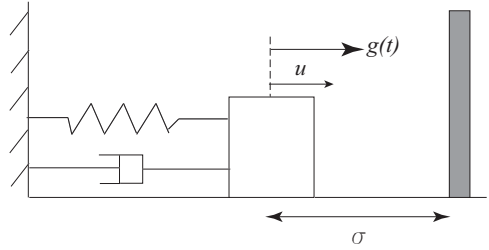}
\caption{A one degree-of-freedom impact oscillator.} \label{fig5}
\end{center}
\end{figure}
\noindent 
The equation of the discontinuity boundary $\Sigma$ in the present case is 
\begin{displaymath}
H({\bf x})=H(u,v)=\sigma-u.
\end{displaymath}
\noindent
Thus we have $\frac{\partial{H}}{\partial{v}}=0$, and hence,
\begin{equation}
\label{eqn8}
H_{\bf x}=(h_1 \;\; 0).
\end{equation}
where $h_1=\frac{\partial{H}}{\partial{u}}$.

The reset map is $R:\Sigma\mapsto\Sigma$, where $R({\bf x})$ has the
form given in (\ref{eqn2}). Note that the
position $u$ does not vary during the impact, i.e., the position of the
mass just before the impact, $u^-$, is same as that just after the
impact, $u^+$, while the velocity of motion $v$ changes. Thus from
(\ref{eqn2})
\begin{equation}
\label{eqn9}
W=(0 \;\; w_2)^T;
\end{equation}
where $w_2$ is the constant of restitution.

\noindent
Using (\ref{eqn8}) and (\ref{eqn9}) in (\ref{eqn6}) we obtain
\begin{eqnarray}
\label{eqn10}
J_{ZDM} & = & {I+{\sqrt{2a^*}}.{\frac{W^*H_{\bf x}}{2\sqrt{-H_{\rm min}}}}} \nonumber \\
 & = & { \left( \begin{array}{cc} 1 & 0 \\ 0 & 1 \end{array} \right) +\frac{\sqrt{2a^*}}{2\sqrt{-H_{\rm min}}}. \left( \begin{array}{c} 0 \\ w_2 \end{array} \right)   \left( \begin{array}{cc} h_1 & 0 \end{array} \right) } \nonumber \\
 & = & \left( \begin{array}{cc} 1 & 0 \\ 0 & 1 \end{array} \right)+\frac{\sqrt{2a^*}}{2\sqrt{-H_{\rm min}}}\left( \begin{array}{cc} 0 & 0 \\ w_2h_1 & 0 \end{array} \right) \nonumber \\
 & = & \left( \begin{array}{cc} 1 & 0 \\ \alpha & 1 \end{array}
\right),
\end{eqnarray}
where \[\alpha = \frac{w_2h_1\sqrt{2a^*}}{2\sqrt{-H_{\rm min}}}\]

\subsection{Investigating the Determinant for Singularity}

From (\ref{eqn5}), the determiant of the normal form map near grazing is
\begin{eqnarray}
|J| & = & |N_2||J_{ZDM}||N_1| \nonumber \\
 & = & |N_2||N_1|  \nonumber
\end{eqnarray}
since $|J_{ZDM}|=1$, from (\ref{eqn10}).

Since the singularity is only in the ZDM, and not in the maps $N_1$
and $N_2$, we conclude that the determinant of the normal form map
does not contain the square-root singularity, and remains
invariant in the immediate neighborhood of the grazing orbit. 

\subsection{Investigating the Trace for Singularity}

To obtain the expression for the trace of the Jacobian $J$ in
(\ref{eqn5}), we need to obtain first the expressions for the maps
$P_1$ and $P_2$.

Let us consider a periodic solution to the equation given in (\ref{eqn7}) as
\begin{displaymath}
p(t)=\left(u(t),v(t)\right)^T.
\end{displaymath}
\noindent
Let $(p(t)+{\delta}p(t))$ be a perturbed orbit, where ${\delta}u$
satisfies the following variational equation
\begin{equation}
\label{eqn11}
{\delta}\ddot{u}+2\zeta{\omega}_n{\delta}\dot{u}+{{\omega}_n}^2{\delta}u=0.
\end{equation}
\noindent
The variational equation needs to be solved to obtain the perturbed
flow ${\delta}p(\tau)=({\delta}u(\tau),{\delta}v(\tau))^T$. Solving
the variational equation amounts to solving the first-order
differential equations
\begin{displaymath}
\frac{d}{dt}{\left( \begin{array}{c} {\delta}u(t)
    \\ {\delta}v(t) \end{array} \right)}= \left( \begin{array}{cc} 0 & 1 \\ -{{\omega}_n}^2 & -2\zeta{{\omega}_n} \end{array}\right) {\left( \begin{array}{c}
    {\delta}u(t) \\ {\delta}v(t) \end{array} \right)}, 
\end{displaymath}
\noindent
with ${\delta}u(0)={\delta}u_0,{\delta}v(0)={\delta}v_0$.

\noindent
The solution of the above problem can be expressed as
\begin{displaymath}
\left( \begin{array}{c} {\delta}u(\tau) \\ {\delta}v(\tau) \end{array} \right)=N_{\tau}\left( \begin{array}{c} {\delta}u_0 \\ {\delta}v_0 \end{array} \right)
\end{displaymath}
\noindent
where 
\begin{equation}
\label{eqn12}
N_{\tau}=e^{-\zeta{\omega}_n\tau} \left( \begin{array}{cc}
  \cos({\omega}_0\tau)+\frac{\zeta}{\sqrt{1-{\zeta}^2}}\sin({\omega}_0\tau)
  & {\sin({\omega}_0\tau)}/{{\omega}_0}
  \\ -\frac{1}{\sqrt{1-{\zeta}^2}}{\omega}_0\sin({\omega}_0\tau) &
  \cos({\omega}_0\tau)-\frac{\zeta}{\sqrt{1-{\zeta}^2}}\sin({\omega}_0\tau) \end{array}\right)
\end{equation}
\noindent
with ${\omega}_0={\omega}_n\sqrt{1-{\zeta}^2}$.

Now we can proceed to obtain the expression of the trace of the
Jacobian in (\ref{eqn5}). In the situation shown in
Fig.~\ref{impact4}, the flow takes time $s_0$ to reach the
discontinuity boundary starting from the Poincar\'e plane, and time
$(T-s_0)$ to return to the Poincar\'e plane starting from the
discontinuity boundary, where $T$ is the time period of the external
forcing function $g(t)$. Thus using the notation in (\ref{eqn12}),
\begin{eqnarray}
N_1 & = & N_{s_0}=e^{-\zeta{\omega}_n{s_0}} \left( \begin{array}{cc} n_{11} & n_{12} \\ n_{13} & n_{14} \end{array}\right) \label{eqn13a} \\
N_2 & = & N_{(T-s_0)}=e^{-\zeta{\omega}_n{(T-s_0)}}
\left( \begin{array}{cc} n_{21} & n_{22} \\ n_{23} &
  n_{24} \end{array} \right) \label{eqn13b}
\end{eqnarray}
\noindent
where the expressions for $n_{11},n_{12},n_{13},n_{14},n_{21},n_{22},n_{23}$ and $n_{24}$ can be derived from (\ref{eqn12}).

\noindent
Using (\ref{eqn10}), (\ref{eqn13a}) and (\ref{eqn13b}) in
(\ref{eqn5}), we get
\begin{eqnarray}
\label{eqn14}
J & = & e^{-\zeta{\omega}_nT}{\left( \begin{array}{cc} n_{21} & n_{22} \\ n_{23} & n_{24} \end{array} \right)}{\left( \begin{array}{cc} 1 & 0 \\ \alpha & 1 \end{array} \right)}{\left( \begin{array}{cc} n_{11} & n_{12} \\ n_{13} & n_{14} \end{array}\right)} \nonumber \\
 & = & e^{-\zeta{\omega}_nT}{\left( \begin{array}{cc} n_{21}+\alpha{n_{22}} & n_{22} \\ n_{23}+\alpha{n_{24}} & n_{24} \end{array} \right)}{\left( \begin{array}{cc} n_{11} & n_{12} \\ n_{13} & n_{14} \end{array} \right)} \nonumber \\
 & = & e^{-\zeta{\omega}_nT}{\left( \begin{array}{cc} n_{21}n_{11}+n_{22}n_{13}+\alpha{n_{22}}{n_{11}} & ^* \\ ^* & n_{23}n_{12}+n_{24}n_{14}+\alpha{n_{24}}n_{12} \end{array} \right)} \nonumber \\
{\Rightarrow}\;Tr(J) & = & e^{-\zeta{\omega}_nT}\{n_{21}n_{11}+n_{22}n_{13}+n_{23}n_{12}+n_{24}n_{14}+{\alpha}\fbox{$(n_{22}n_{11}+n_{24}n_{12})$}\}
\end{eqnarray}

The expression for the trace of the Jacobian $Tr(J)$, in
(\ref{eqn14}), shows that the singularity term $\alpha$ has a
coefficient $e^{-\zeta{\omega}_nT}(n_{22}n_{11}+n_{24}n_{12})$. Let us
take a closer look at this coefficient. Using (\ref{eqn12}),
(\ref{eqn13a}) and (\ref{eqn13b})
\begin{eqnarray}
n_{11} & = & \cos({\omega}_0s_0)+\frac{\zeta}{\sqrt{1-{\zeta}^2}}\sin({\omega}_0s_0) \nonumber \\
n_{12} & = & \frac{\sin({\omega}_0s_0)}{{\omega}_0} \nonumber \\
n_{22} & = & \frac{\sin\{{\omega}_0(T-s_0)\}}{{\omega}_0} \nonumber \\
and{\:}n_{24} & = & \cos\{{\omega}_0(T-s_0)\}-\frac{\zeta}{\sqrt{1-{\zeta}^2}}\sin\{{\omega}_0(T-s_0)\}. \nonumber
\end{eqnarray}
\noindent
Therefore it follows that
\begin{equation}
\label{eqn15}
\fbox{$n_{22}n_{11}+n_{24}n_{12}=\frac{\sin({\omega}_0T)}{{\omega}_0}\neq{0}$},{\;\;}
\forall{\;}{\omega}_0\neq{\frac{m{\omega}_{\rm forcing}}{2}}
\end{equation}
\noindent
where ${\omega}_{\rm forcing}$ is the angular frequency of the
periodic forcing function $g(t)$, i.e., ${\omega}_{\rm
  forcing}T=2\pi$, and $m\in{I}$.  Thus the coeffcient of $\alpha$ in
the expression of the trace (\ref{eqn14}) of the Jacobian of the
stroboscopic map must be a non-zero entity. Thus the singularity in
$\alpha$ survives, and hence a square-root singularity must occur in
the trace of the Jacobian.

\section{Conclusions}

In this paper we have probed the variation of the trace and
determinant of the Jacobian matrix of map of a hard-impact oscillator
as it goes from non-impacting state to an impacting state. We have
shown that the square-root singularity should be expressed only in the
trace of the Jacobian matrix while the determinant should remain
invariant in the immediate neighborhood of a grazing orbit.


\end{document}